\documentclass{elsart3-1}

\usepackage{amssymb}
\usepackage{amsmath,latexsym,graphicx}
\usepackage[english,francais]{babel}

\newtheorem{theorem}{Theorem}[section]
\newtheorem{lemma}[theorem]{Lemma}
\newtheorem{e-proposition}[theorem]{Proposition}
\newtheorem{corollary}[theorem]{Corollary}
\newtheorem{e-definition}[theorem]{Definition\rm}

\newtheorem{theoreme}{Th\'eor\`eme}[section]

\newtheorem{proposition}[theoreme]{Proposition}

\renewcommand{\theequation}{\arabic{equation}}
\newcounter{count}

\DeclareMathOperator{\sech}{sech}

\def\og{\leavevmode\raise.3ex\hbox{$\scriptscriptstyle\langle\!\langle$~}}
\def\fg{\leavevmode\raise.3ex\hbox{~$\!\scriptscriptstyle\,\rangle\!\rangle$}}

\makeatletter
\newcommand{\defeq}{\mathrel{\rlap{\raisebox{0.3ex}{$\m@th\cdot$}}\raisebox{-0.3ex}{$\m@th\cdot$}}=}
\makeatother

\newcommand{\dx}{\mathrm{d}x}

\newcommand{\R}{\mathbb R}

\newcommand{\Z}{\mathbb Z}

\renewcommand{\AA}{\mathcal A}
\newcommand{\BB}{\mathcal B}
\newcommand{\CC}{\mathcal C}

\newcommand{\LL}{\mathcal L}
\newcommand{\UU}{\mathcal U}

\newcommand{\XX}{\mathcal X}
\newcommand{\YY}{\mathcal Y}
\newcommand{\ZZ}{\mathcal Z}

\renewcommand{\e}{\mathrm{e}}
\renewcommand{\i}{\mathrm{i}}

\newcommand{\eqn}[1]{(\ref{#1})}
\journal{the Acad\'emie des sciences}
\begin{document}
% place in the next line the header (rubrique) chosen for your article,
% if you know it (you can also have 2, format : Header1/Header2
\centerline{}
\vspace{-4cm}
\begin{frontmatter}


% \author{Name\thanksref{label2}}
\selectlanguage{english}
\title{Periodic solitons for the elliptic-elliptic\\ focussing Davey-Stewartson equations}

% The [label1] can be suppressed if there is only one address for all authors

\selectlanguage{english}
\author[sb,leic]{Mark D. Groves},
\ead{groves@math.uni-sb.de}
\author[vt]{Shu-Ming Sun},
\ead{sun@math.vt.edu}
\author[lund]{Erik Wahl\'{e}n}
\ead{ewahlen@maths.lth.se}

\address[sb]{FR 6.1 - Mathematik, Universit\"{a}t des Saarlandes, Postfach 151150, 66041 Saarbr\"{u}cken, Germany}
\address[leic]{Department of Mathematical Sciences, Loughborough University, Loughborough, Leics, LE11 3TU, UK}
\address[vt]{Department of Mathematics, Virginia Polytechnic Institute and
State University, Blacksburg, VA 24061, USA}
\address[lund]{Department of Mathematics, Lund University, 22100 Lund, Sweden}

% If you know the dates of reception, and acceptation you can put them now;
%  idem the name of the person presenting the Note

\medskip
\begin{center}
{\small Received *****; accepted after revision +++++\\
Presented by *****}
\end{center}

\begin{abstract}
\selectlanguage{english}
We consider the elliptic-elliptic, focussing Davey-Stewartson equations, which have an explicit
bright line soliton solution. The existence of a
family of periodic solitons, which have the profile of the line soliton in the longitudinal spatial
direction and are periodic in the transverse spatial direction, is established using
dynamical systems arguments. We also show that the line soliton is linearly unstable
with respect to perturbations in the transverse direction.

\vskip 0.5\baselineskip

\selectlanguage{francais}
% Text of abstract in French
\noindent{\bf R\'esum\'e} \vskip 0.5\baselineskip \noindent
{\bf Solitons p\'{e}riodiques du syst\`{e}me de Davey-Stewartson elliptique-elliptique focalisant.}
Nous consid\'{e}rons les \'{e}quations de Davey-Stewartson focalisantes dans le cas elliptique-elliptique, lorsqu'elles poss\`{e}dent une solution unidimensionnelle de type soliton. En utilisant des m\'{e}thodes de la th\'{e}orie des syst\`{e}mes dynamiques, nous montrons l'existence d'une famille de solutions bidimensionnelles qui ont le profil d'un soliton dans la direction spatiale longitudinale et sont p\'{e}riodiques dans la direction spatiale transverse. Nous montrons \'{e}galement que le soliton unidimensionnel est lin\'{e}airement instable vis-\`{a}-vis des perturbations transverses.

\end{abstract}
\end{frontmatter}

\selectlanguage{english}
% main text

\vspace*{-1.3cm}

\section{Introduction}

The Davey-Stewartson equations
\begin{align}
\i A_t + \varepsilon A_{xx} +  A_{yy} +(\gamma_1 |A|^2 +\gamma_2 \phi_x)A&=0, \label{DS1} \\
\gamma_3\phi_{xx} + \phi_{yy} -\gamma_3 |A|_x^2 &= 0, \label{DS2}
\end{align}
where $\varepsilon = \pm 1$, $\gamma_1$, $\gamma_2$, $\gamma_3 \in {\mathbb R} \setminus \{0\}$
with $\gamma_2+\gamma_3 = \pm 2$, arise in the modelling of wave packets on the surface
of a three-dimensional body of water; the variables $A=A(x,y,t)$ and
$\phi=\phi(x,y,t)$ are the complex wave amplitude and real mean flow and the signs of the
parameters depend upon the particular physical regime under consideration (see
Ablowitz \& Segur \cite[\S2.2]{AblowitzSegur79}). In the literature the cases $\gamma_1+\gamma_2=2$ and
$\gamma_1+\gamma_2=-2$ are termed respectively \emph{focussing} and \emph{defocussing}, and the
system is classified as hyperbolic-hyperbolic, hyperbolic-elliptic, elliptic-hyperbolic or
elliptic-elliptic according to the signs of $\varepsilon$ and $\gamma_3$. Certain special cases
of the mixed-type systems are often referred to as DS-I and DS-II
and are known to be completely integrable (see Ablowitz \& Clarkson \cite[p.\ 60]{AblowitzClarkson}). Note that
solutions of \eqn{DS1}, \eqn{DS2} which are spatially homogeneous in the $y$-direction
satisfy the cubic nonlinear Schr\"{o}dinger equation
\begin{equation}
\i A_t + A_{xx} +(\gamma_1+\gamma_2)|A|^2A = 0 \label{NLS}
\end{equation}
(where $\phi$ is recovered from \eqn{DS2}).\pagebreak

Solutions of \eqn{NLS} which converge to an equilibrium as $x \rightarrow \pm \infty$
and are $2\pi$-periodic in $t$ are referred to as \emph{line solitons}. In the defocussing
case the equation admits a `dark' line soliton which decays to a nontrivial equilibrium,
while in the focussing case it has the `bright' line soliton
\begin{equation}
A^\star(x,t) = \e^{\i t} \sech(x) \label{BLS}
\end{equation}
which satisfies $A^\star(x,t) \rightarrow 0$ as $x \rightarrow \pm\infty$.
In this note we examine \emph{periodic solitons} which decay as
$x \rightarrow \pm \infty$ and are periodic in $y$ and $t$, and in particular consider how
they emerge from line solitons in a \emph{dimension-breaking bifurcation}. Explicit
formulae for dark periodic solitons have been obtained for the integrable versions of the
equations by Watanabe \& Tajiri \cite{WatanabeTajiri98} and Arai, Takeuchi \& Tajiri \cite{AraiTakeuchiTajiri01}; here we establish
the existence of bright periodic solitons to the elliptic-elliptic, focussing
equations ($\varepsilon=1$, $\gamma_1+\gamma_2=2$, $\gamma_3>0$) under the additional condition
$\gamma_2>0$ by dynamical-systems methods.

\begin{theorem} \label{Main result 1}
Suppose that $\varepsilon=1$, $\gamma_1+\gamma_2=2$ and $\gamma_2$, $\gamma_3>0$.
There exist an open neighbourhood
${\mathcal N}$ of the origin in ${\mathbb R}$, a positive real number $\omega_0$
and a family of periodic solitons
$\{\e^{\i t}u_s(x,y), \phi_s(x,y)
\}_{s \in {\mathcal N}}$ to \eqn{DS1}, \eqn{DS2} which
emerges from the bright line soliton in
a dimension-breaking bifurcation. Here
\[u_s(x,y)=\sech(x) + u_s^\prime(x,y), \qquad \phi_s(x,y) = \tanh(x) + \phi_s^\prime(x,y),\]
in which $u_s^\prime(\cdot\,,\cdot)$, $\phi_s^\prime(\cdot\,,\cdot)$ are real, have amplitude $O(|s|)$
and are even in both arguments and periodic in their second with frequency $\omega_0 + O(|s|^2)$.
\end{theorem}

We also present a corollary to this result which asserts that the bright line
soliton is \emph{transversely linearly unstable} and thus confirms the prediction
made by Ablowitz \& Segur \cite[\S 3.2]{AblowitzSegur79}.

\begin{theorem} \label{Main result 2}
Suppose that $\varepsilon=1$, $\gamma_1+\gamma_2=2$ and $\gamma_2$, $\gamma_3>0$. For
each sufficiently small positive value of $\lambda$ the linearisation of \eqn{DS1}, \eqn{DS2} at
$A^\star(x,t) = \e^{\i t} \sech(x)$, $\phi^\star(x,y)=\tanh(x)$ has a solution of the form
$\e^{\lambda t+\i t}(A(x,y),\phi(x,y))$, where $(A(x,y),\phi(x,y))$ is periodic in $y$ and satisfies
$(A(x,y),\phi(x,y))\rightarrow (0,0)$ as $x \rightarrow \pm \infty$.
\end{theorem}

In the remainder of this article we suppose that $\varepsilon=1$, $\gamma_1+\gamma_2=2$ and $\gamma_2$, $\gamma_3>0$.
Equations \eqn{DS1}, \eqn{DS2} with these coefficients arise when modelling
water waves with weak surface tension. The existence and transverse linear instability of periodic solitons for the water-wave problem in this physical regime
has recently been established by Groves, Sun \& Wahl\'{e}n
\cite{GrovesSunWahlen16}.

\vspace*{-1cm}

\section{Spatial dynamics}

The equations for solutions of \eqn{DS1}, \eqn{DS2} for which $A(x,y,t)=\e^{\i t} \big(u_1(x,y,t)+ \i u_2(x,y,t)\big)$
(and  $u_1$, $u_2$ are real-valued) can be formulated as the evolutionary system
\begin{align}
u_{1y} & = v_1, \label{Full SD1} \\
v_{1y} & =  u_{2t}  - u_{1xx}  +u_1 -(\gamma_1 u_1^2 + \gamma_1 u_2^2 +\gamma_2 \phi_x)u_1, \label{Full SD2} \\
u_{2y} & = v_2, \label{Full SD3} \\
v_{2y} & = -u_{1t} - u_{2xx}  +u_2 -(\gamma_1 u_1^2 + \gamma_1 u_2^2 +\gamma_2 \phi_x)u_2, \label{Full SD4} \\
\phi_y & = \psi, \label{Full SD5} \\
\psi_y & = - \gamma_3\phi_{xx} +\gamma_3 (u_1^2+u_2^2)_x, \label{Full SD6}
\end{align}
where the spatial direction $y$ plays the role of time. 
To identify an appropriate functional-analytic setting for these equations,
let us first specialise to stationary solutions, so that
\begin{align}
u_{1y} & = v_1, \label{SD1} \\
v_{1y} & =  - u_{1xx}  +u_1 -(\gamma_1 u_1^2 + \gamma_1 u_2^2 +\gamma_2 \phi_x)u_1, \label{SD2} \\
u_{2y} & = v_2, \label{SD3} \\
v_{2y} & = - u_{2xx}  +u_2 -(\gamma_1 u_1^2 + \gamma_1 u_2^2 +\gamma_2 \phi_x)u_2, \label{SD4} \\
\phi_y & = \psi, \label{SD5} \\
\psi_y & = - \gamma_3\phi_{xx} +\gamma_3 (u_1^2+u_2^2)_x. \label{SD6}
\end{align}

Equations \eqn{SD1}--\eqn{SD6}
constitute a semilinear evolutionary system in the phase space 
$X=H^1(\R) \times L^2(\R) \times H^1(\R) \times L^2(\R) \times H^1(\R) \times L^2(\R)$; the domain of the
linear part of the vector field defined by their right-hand side is
$D=H^2(\R) \times H^1(\R) \times H^2(\R) \times H^1(\R) \times H^2(\R) \times H^1(\R)$.
This evolutionary system is
reversible, that is invariant under $y \mapsto -y$,
$(u_1,v_1,u_2,v_2,\phi,\psi) \mapsto
S(u_1,v_1,u_2,v_2,\phi,\psi)$, where the reverser $S: X \rightarrow X$ is defined by
$S(u_1,v_1,u_2,v_2,\phi,\psi)=(u_1,-v_1,u_2,-v_2,\phi,-\psi)$. It is also invariant under the reflection
$R: X \rightarrow X$ given by $R(u_1(x),v_1(x),u_2(x),v_2(x),\phi(x),\psi(x))=(u_1(-x),v_1(-x),u_2(-x),v_2(-x),-\phi(-x),-\psi(-x))$,
and one may seek solutions which are invariant under this symmetry by replacing $X$
and $D$ by respectively
$$
X_\mathrm{r}:=X \cap \mathrm{Fix}\, R = H_\mathrm{e}^1(\R) \times L_\mathrm{e}^2(\R) \times H_\mathrm{e}^1(\R) \times L_\mathrm{e}^2(\R) \times H_\mathrm{o}^1(\R) \times L_\mathrm{o}^2(\R)
$$
and
$$
D_\mathrm{r}:=D \cap \mathrm{Fix}\, R = H_\mathrm{e}^2(\R) \times H_\mathrm{e}^1(\R) \times H_\mathrm{e}^2(\R) \times H_\mathrm{e}^1(\R) \times H_\mathrm{o}^2(\R) \times H_\mathrm{o}^1(\R),
$$
where
\begin{align*}
H_\mathrm{e}^n(\R) & = \{w \in H^n(\R): w(x)=w(-x) \text{ for all }  x \in \R\}, \\
H_\mathrm{o}^n(\R) & = \{w \in H^n(\R): w(x)=-w(-x) \text{ for all }  x \in \R\}.
\end{align*}
It is also possible to replace $D_\mathrm{r}$ by the extended function space
$$
D_\star:= H_\mathrm{e}^2(\R) \times H_\mathrm{e}^1(\R) \times H_\mathrm{e}^2(\R) \times H_\mathrm{e}^1(\R) \times H_{\star,\mathrm{o}}^2(\R) \times H_\mathrm{o}^1(\R),
$$
where
$$H_{\star,\mathrm{o}}^2(\R)=\{w\in L_\mathrm{loc}^2(\R): w_x \in H^1(\R), w(x)=-w(-x) \text{ for all }  x \in \R\}
$$
(a Banach space with norm
$\|w\|_{\star, 2}:=\|w_x\|_1$).
This feature allows one to consider solutions to \eqn{SD1}--\eqn{SD6} whose $\phi$-component is not evanescent;
in particular solutions corresponding to line solitons fall into this category (see below).

\begin{figure}[h]
\begin{center}
\includegraphics[width=9cm]{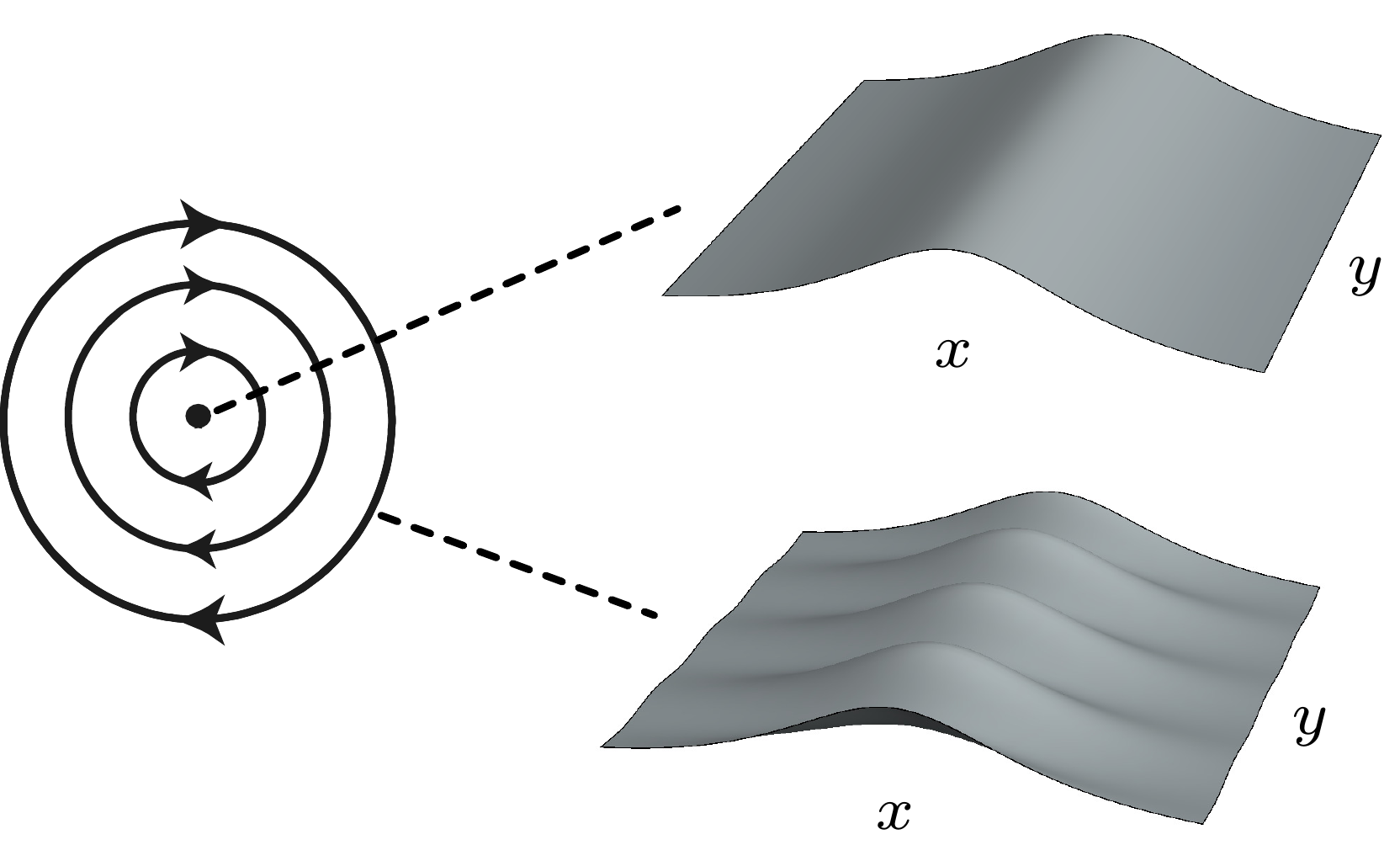}
\end{center}
\caption{A family of periodic solutions surrounding a nontrivial equilibrium solution to \eqn{SD1}--\eqn{SD6}
in its phase space (left) corresponds to a dimension-breaking bifurcation of a branch of periodic solitons from
a line soliton (right, plot of $u(x,y)$).} \label{dbreak}
\end{figure}

Each point in phase space
corresponds to a function on the real line which decays as $x \rightarrow \infty$, and the dynamics
of equations \eqn{SD1}--\eqn{SD6} in $y$
describes the behaviour of their solutions in the $y$-direction. In particular, equilibria
correspond to line solitons (the equilibrium
$$(u_1^\star(x),v_1^\star(x),u_2^\star(x),v_2^\star(x),\phi^\star(x),\psi^\star(x))=(\sech(x), 0,0,0,\tanh(x), 0)$$
corresponds to the line soliton \eqn{BLS}),
while periodic orbits correspond to periodic solitons (see Figure \ref{dbreak}).
In Section \ref{Apply LI} we construct dimension-breaking bifurcations by writing
\begin{equation}
(u_1,v_1,u_2,v_2,\phi,\psi)=(u_1^\star,v_1^\star,\phi^\star,u_2^\star,v_2^\star,\phi^\star,\psi^\star)+(u_1^\prime,v_1^\prime,u_2^\prime,v_2^\prime,\phi^\prime,\psi^\prime)
\label{Main Ansatz}
\end{equation}
and seeking small-amplitude periodic solutions of the resulting evolutionary system
\begin{equation}
w_y=Lw + N(w) \label{DS final}
\end{equation}
for $w=(u_1^\prime,v_1^\prime,u_2^\prime,v_2^\prime,\phi^\prime,\psi^\prime)$, where
$$
L\begin{pmatrix}
u_1 \\
v_1 \\
u_2 \\
v_2 \\
\phi \\
\psi
\end{pmatrix}
=
\begin{pmatrix}
v_1 \\
- u_{1xx} + u_1 -(3\gamma_1+\gamma_2) \sech^2 (x) u_1 -\gamma_2 \sech (x)\phi_x \\
v_2 \\
- u_{2xx} + u_2 -2\sech^2 (x) u_2 \\
\psi \\
- \gamma_3\phi_{xx} +2\gamma_3 (\sech (x) u_1)_x
\end{pmatrix}, 
$$
$$
N\begin{pmatrix}
u_1 \\
v_1 \\
u_2 \\
v_2 \\
\phi \\
\psi
\end{pmatrix}
=
\begin{pmatrix}
0 \\
-3\gamma_1\sech (x) u_1^2
-\gamma_1 \sech(x)^2 u_2^2
-\gamma_2u_1\phi_x
-\gamma_1 u_1^3
-\gamma_1 u_1u_2^2\\
0 \\
-2\gamma_1\sech (x) u_1u_2
-\gamma_2u_2\phi_x
-\gamma_1 u_2^3
-\gamma_1 u_1^2u_2\\
0\\
\gamma_3(u_1^2+u_2^2)_x
\end{pmatrix}
$$
and we have dropped the primes for notational simplicity. Note that \eqn{DS final}
has the invariant subspace $\tilde{X}=\{(u_2,v_2)=(0,0)\}$, and we define
$\tilde{X}_\mathrm{r}=X_\mathrm{r} \cap \tilde{X}$, $\tilde{D}_\mathrm{r} = D_\mathrm{r}
\cap \tilde{X}$ and $\tilde{D}_\star = D_\star \cap \tilde{X}$.

Returning to \eqn{Full SD1}--\eqn{Full SD6}, observe that these equations
constitute a reversible evolutionary equation with phase space $H^1((-t_0,t_0),X)$;
the domain of its vector field is
$H^2((-t_0,t_0),X) \cap H^1((-t_0,t_0), D)$
and its reverser is given by the pointwise extension of $S: X \rightarrow X$ to $H^1((-t_0,t_0),X)$.
Seeking solutions of the form \eqn{Main Ansatz}, we find that
$$w_y = T w_t + Lw+N(w),$$
where $T(u_1,v_1,u_2,v_2,\phi,\psi) = (0, u_2,0,-u_1,0,0)$
and we have again dropped the primes. In Section \ref{Apply GH} we demonstrate that
the solution $(u_1^\star,v_1^\star,u_2^\star,v_2^\star,\phi^\star,\psi^\star)$ of
\eqn{Full SD1}--\eqn{Full SD6} is transversely linearly unstable by constructing
a solution of the linear equation
\begin{equation}
w_y = T w_t + Lw \label{TS}
\end{equation}
of the form $\e^{\lambda t} u_\lambda(y)$, where
$u_\lambda \in C_\mathrm{b}^1({\mathbb R}, X) \cap C_\mathrm{b}(\R, D)$ is
periodic, for each sufficiently small positive value of $\lambda$.

\vspace*{-1cm}

\section{Spectral theory}

In this section we determine the purely imaginary spectrum of the linear operator
$L: D \subseteq X \to X$. To this end we study the resolvent equations
\begin{equation}
(L-\i kI )w= w^\dag \label{Resolvent equations}
\end{equation}
for $L$, where $w=(u_1,v_1,u_2,v_2,\phi,\psi)$, $w^\dag=(u_1^\dag,v_1^\dag,u_2^\dag,v_2^\dag,\phi^\dag,\psi^\dag)$ and
$k \in {\mathbb R} \setminus \{0\}$; since $L$ is real and anticommutes with the reverser $S$ it suffices to examine
non-negative values of $k$, real values of $u_1$, $u_2$, $\phi$, $v_1^\dag$, $v_2^\dag$, $\psi^\dag$ and purely
imaginary values of $u_1^\dag$, $u_2^\dag$, $\phi^\dag$, $v_1$, $v_2$, $\psi$. Observe that \eqn{Resolvent equations}
is equivalent to the decoupled equations
$$(\AA_1+k^2 I)
\begin{pmatrix}
u_1 \\
\phi
\end{pmatrix}
=
\begin{pmatrix}
v_1^\dag + \i k u_1^\dag \\
\psi^\dag + \i k \phi^\dag
\end{pmatrix},
\qquad
(\AA_2+k^2 I) u_2 = v_2^\dag + \i k u_2^\dag,
$$
where $\AA_1: H^2(\R) \times H^2(\R) 
\subseteq  L^2(\R) \times L^2(\R)
\to L^2(\R) \times L^2(\R)$ and $\AA_2: H^2(\R)
\subseteq  L^2(\R) \to L^2(\R)$ are
defined by
\begin{align*}
\AA_1 \begin{pmatrix}
u_1 \\
\phi
\end{pmatrix}
& =
\begin{pmatrix}
-u_{1xx} + u_1 -(3\gamma_1+\gamma_2) \sech^2 (x)u_1 -\gamma_2\sech (x)\phi_x \\
- \gamma_3\phi_{xx} +2\gamma_3 (\sech (x) u_1)_x
\end{pmatrix}, \\[2mm]
\AA_2 u_2 & = -u_{2xx}+u_2 - 2 \sech^2(x)u_2;
\end{align*}
the values of $v_1$, $v_2$ and  $\psi$ are recovered from the formulae
$$v_1=u_1^\dag+\i k u_1, \qquad v_2=u_2^\dag+\i k u_2, \qquad \psi=\phi^\dag+\i k \phi.$$
It follows that $L-\i k I$ is (semi-)Fredholm if $\AA_1+k^2 I$ and $\AA_2+k^2 I$ are (semi-)Fredholm
and the dimension of the (generalised) kernel of $L-\i  k I$ is the sum of those of $\AA_1+k^2 I$ and $\AA_2+k^2 I$.

Lemmata \ref{Spectrum of AA1} and \ref{Spectrum of AA2} below record the spectra of $\AA_1$ and $\AA_2$;
part (i) of the following proposition (see Drazin \cite[Chapter 4.11]{Drazin}) is used in the proof of the former
while the latter follows directly from part (ii).

\begin{proposition} \label{Drazin prop}\hspace{1cm}
\begin{list}{(\roman{count})}{\usecounter{count}}
\item
The spectrum of the self-adjoint operator
$1-\partial_x^2-6 \sech^2(x): H^2(\R)\subseteq L^2(\R)\to L^2(\R)$ 
consists of essential spectrum $[1,\infty)$ and two simple eigenvalues at $-3$ and $0$ (with corresponding eigenvectors
$\sech^2(x)$ and $\sech^\prime(x)$).
\item
The spectrum of the self-adjoint operator
$1-\partial_x^2-2\sech^2(x): H^2(\R)\subseteq L^2(\R)\to L^2(\R)$ 
consists of essential spectrum $[1,\infty)$ and a simple eigenvalue at $0$ (with corresponding eigenvector $\sech(x)$).
\end{list}
\end{proposition}

\begin{lemma} \label{Spectrum of AA1}
The spectrum of the operator $\AA_1$ consists of essential spectrum $[0,\infty)$ and an algebraically simple negative eigenvalue
$-\omega_0^2$ whose eigenspace lies in $L_\mathrm{e}^2(\R) \times L_\mathrm{o}^2(\R)$.
\end{lemma}
{\bf Proof.} First note that $\AA_1$ is a compact perturbation of the constant-coefficient operator\linebreak
$H^2(\R) \times H^2(\R) 
\subseteq L^2(\R) \times L^2(\R)
\to L^2(\R) \times L^2(\R)$
defined by
\[
(u_1,\phi)\mapsto (-u_{1xx} + u_1, -\gamma_3\phi_{xx}),
\]
whose essential spectrum is clearly $[0,\infty)$; it follows that $\sigma_\mathrm{ess}(\AA_1) = [0,\infty)$
(see Kato, \cite[Chapter IV, Theorem 5.26]{Kato}). Because $\AA_1$ is self-adjoint with respect to the inner product $\langle (u_1^1,\phi^1), (u_1^2,\phi^2) \rangle
= \langle u_1^1, u_1^2 \rangle_0 + \frac{1}{2}\gamma_2\gamma_3^{-1} \langle \phi^1, \phi^2 \rangle_0$ for $L^2(\R) \times L^2(\R)$
the remainder of its spectrum
consists of negative real eigenvalues with finite multiplicity. 

One finds by
an explicit calculation that
$$
\langle \AA_1(u_1,\phi), (u_1,\phi)\rangle
=
\langle u_1-u_{1xx}-6 \sech^2 (x)u_1,u_1 \rangle_0
+\frac{\gamma_2}{2}\int_\R\left(\phi_x-2\sech (x) u_1\right)^2\dx,
$$
which quantity is positive for $(u_1,\phi) \in H^2(\R) \times H^2(\R) $ with
$$\langle (u_1,\phi), (\sech^2(x),0)\rangle=0$$
(see Proposition \ref{Drazin prop}(ii)).
It follows that any subspace of $H^2(\R) \times H^2(\R) $ upon which $\AA_1$
is strictly negative definite is one-dimensional.
The calculation
\[
\lim_{R \rightarrow \infty} \langle \AA_1 (\sech(x), 2\phi_R(x)\tanh(x)),(\sech(x),2\phi_R(x)\tanh(x))\rangle=-\frac{16}{3},
\]
where $\phi(R)=\chi(x/R)$ and $\chi \in C_0^\infty({\mathbb R})$ is a smooth cut-off function equal to
unity in $[-1,1]$, shows that $\inf \sigma(\AA_1)<0$, so that the spectral subspace of $H^2(\R) \times H^2(\R) $
corresponding to the part of the spectrum of $\AA_1$ in $(-\infty,-\varepsilon)$
is nontrivial and hence one-dimensional for every sufficiently small value of $\varepsilon>0$.
We conclude that $\AA_1$ has precisely one simple negative eigenvalue
$-\omega_0^2$.

Finally, the same argument shows that $\AA_1|_{L_\mathrm{e}^2(\R) \times L_\mathrm{o}^2(\R)}$
also has precisely one simple negative eigen- value. It follows that this eigenvalue is $-\omega_0^2$,
whose eigenspace therefore lies in $L_\mathrm{e}^2(\R) \times L_\mathrm{o}^2(\R)$.\qed

\begin{lemma} \label{Spectrum of AA2}
The spectrum of the operator $\AA_2$ consists of essential spectrum $[1,\infty)$ and an algebraically simple negative eigenvalue
at $0$ whose eigenspace lies in $L_\mathrm{e}^2(\R)$.
\end{lemma}

\begin{corollary} \label{Imaginary result}
The purely imaginary number $\i k$ belongs to the resolvent set of $L$ for $k \in \R \setminus \{0, \pm \omega_0\}$
and $\pm \i \omega_0$ are algebraically simple purely imaginary eigenvalues of $L$ whose eigenspace lies in
$\tilde{X}_\mathrm{r}$.
\end{corollary}

\section{Application of the Lyapunov-Iooss theorem} \label{Apply LI}

Our existence theory for periodic solitons
is based upon an application of the following version
of the Lyapunov centre theorem for reversible systems (see Iooss \cite{Iooss99})
which allows
for a violation of the classical nonresonance condition at the origin
due to the presence of essential spectrum there (a feature typical of spatial
dynamics formulations for problems in unbounded domains)
provided that the `Iooss condition at the origin' (hypothesis (viii)) is satisfied.

\begin{theorem}[Iooss-Lyapunov]
Consider the differential equation
\begin{equation}
w_\tau=L(w)+N(w), \label{Basic DE}
\end{equation}
in which $w(\tau)$ belongs to a real Banach space $\XX$. Suppose that $\YY$, $\ZZ$ are further real Banach spaces
with the properties that
\begin{list}{(\roman{count})}{\usecounter{count}}
\item
$\ZZ$ is continuously embedded in $\YY$ and continuously and densely embedded in $\XX$,
\item
$L: \ZZ \subseteq \XX \rightarrow \XX$ is a closed linear operator,
\item
there is an open neighbourhood $\UU$ of the origin in $\YY$
such that $L \in \LL(\YY,\XX)$ and $N \in C^3_{\mathrm{b,u}}(\UU, \XX)$ (and hence $N \in C^3_{\mathrm{b,u}}(\UU \cap \ZZ, \XX)$)
with $N(0)=0$, $\mathrm{d}N[0]=0$.
\end{list}
Suppose further that
\begin{list}{(\roman{count})}{\usecounter{count}}
\setcounter{count}{3}
\item equation \eqn{Basic DE} is reversible: there exists an involution $S \in \LL(\XX)
\cap \LL(\YY) \cap \LL(\ZZ)$
with $SLw=-LSw$ and
$SN(w)=-N(Sw)$ for all $w \in \UU$,
\end{list}
and that the following spectral hypotheses are satisfied.
\begin{list}{(\roman{count})}{\usecounter{count}}
\setcounter{count}{4}
\item $\pm \i \omega_0$ are nonzero simple eigenvalues of $L$;
\item $\i n \omega_0 \in \rho(L)$ for $n \in \Z\!\setminus\!\{-1,0,1\}$;
\item $\|(L-\i n \omega_0 I )^{-1}\|_{\XX \rightarrow \XX} = o(1)$ and $\|(L-\i n \omega_0 I )^{-1}\|_{\XX \rightarrow \ZZ} = O(1)$ as $n \rightarrow \pm\infty$;
\item For each $w^\dag \in \UU$ the equation
$$Lw=-N(w^\dag)$$
has a unique solution $w \in \YY$ and the mapping $w^\dag \mapsto w$ belongs to $C^3_{\mathrm{b,u}}(\UU, \YY)$.
\end{list}
Under these hypotheses there exist an open neighbourhood $I$ of the origin in $\R$ and a\linebreak
continuously differentiable branch $\{(v(s),\omega(s))\}_{s \in I}$ of reversible, $2\pi/\omega(s)$-periodic solutions in\linebreak
$C^1_\mathrm{per}({\mathbb R}, \YY \oplus \XX) \cap C_\mathrm{per}({\mathbb R}, \YY \oplus \ZZ)$
to \eqn{Basic DE} with amplitude $O(|s|)$.
Here the direct sum refers to the decomposition of a function into its mode $0$ and higher-order Fourier components,
the subscript `per' indicates a $2\pi/\omega(s)$-periodic function and $\omega(s)=\omega_0+O(|s|^2)$.
\end{theorem}

Theorem \ref{Main result 1} is proved by applying the Iooss-Lyapunov theorem
to \eqn{DS final}, taking $\XX=\tilde{X}_\mathrm{r}$, $\YY=\tilde{D}_\star$, $\ZZ=\tilde{D}_\mathrm{r}$ and $\UU=\tilde{D}_\star$
(and of course $\tau=y$ and $S(u_1,v_1,\phi,\psi)=(u_1,-v_1,\phi,-\psi)$).
The spectral hypotheses (v) and (vi) follow from Corollary \ref{Imaginary result}, while
(vi) and (vii) are verified in respectively Lemma \ref{Resolvent estimate} and \ref{Iooss condition}
below.

\begin{lemma} \label{Resolvent estimate}
The operator $L$ satisfies the resolvent estimates 
$\|(L-\i k I )^{-1}\|_{\tilde{X}_\mathrm{r} \rightarrow \tilde{X}_\mathrm{r}} = O(|k|^{-1})$ and\linebreak
$\|(L-\i k I )^{-1}\|_{\tilde{X}_\mathrm{r} \rightarrow \tilde{D}_\mathrm{r}} = O(1)$ as $|k| \rightarrow \infty$.
\end{lemma}
{\bf Proof.} Notice that
$$L(u_1,v_1,\phi,\psi)=(\BB_1(u_1,v_1),\BB_2(\phi,\psi))+\CC(u_1,v_1,\phi,\psi),$$
where $\BB_1(u_1,v_1)=(v_1,-u_{1xx}+u_1)$, $\BB_2(\phi,\psi)=(\psi,-\gamma_3\phi_{xx}+\gamma_3\phi)$
and
$$\CC(u_1,v_1,\phi,\psi) = (0,-(3\gamma_1+\gamma_2) \sech^2 (x) u_1 - \gamma_2\sech (x)\phi_x,0, -\gamma_3\phi +2\gamma_3 (\sech (x) u_1)_x).$$

Writing $\tilde{X}_\mathrm{r} = \tilde{X}_1 \times \tilde{X}_2$, $\tilde{D}_r = \tilde{Y}_1 \times \tilde{Y}_2$ and equipping $\tilde{X}_1$ with the usual inner product
and $\tilde{X}_2$ with the inner product $\langle (\phi^1,\psi^1), (\phi^2,\psi^2) \rangle
= \langle \phi^1, \phi^2 \rangle_1 + \gamma_3^{-1} \langle \psi^1, \psi^2 \rangle_0$, observe that
$\BB_j: \tilde{Y}_j \subset \tilde{X}_j \to \tilde{X}_j$ is self-adjoint, so that
$$\|(\BB_j-\i k I)^{-1}\|_{\tilde{X}_j \to \tilde{X}_j} \leq |k|^{-1}$$
for $k \neq 0$.
Furthermore
$\|\BB_j (\cdot)\|_{\tilde{X}_j} = \|\cdot\|_{\tilde{Y}_j}$,
so that
$$
\|(\BB_j-\i k I)^{-1}\|_{\tilde{X}_j \to \tilde{Y}_j} = \|\BB_j (\BB_j-\i k I)^{-1}\|_{\tilde{X}_j \to \tilde{X}_j}=\| I + \i k I(\BB_j-\i k I)^{-1}\|_{\tilde{X}_j \to \tilde{X}_j} \leq 2
$$
for $k \neq 0$.
It follows that $\BB=\BB_1 \times \BB_2: \tilde{D}_\mathrm{r} \subseteq \tilde{X}_\mathrm{r} \to \tilde{X}_\mathrm{r}$ satisfies the estimates
\begin{equation}
\|(\BB-\i k I)^{-1}\|_{\tilde{X}_\mathrm{r} \to \tilde{X}_\mathrm{r}} \leq |k|^{-1}, \qquad \|(\BB-\i k I)^{-1}\|_{\tilde{X}_\mathrm{r} \to \tilde{D}_\mathrm{r}} \leq 2
\label{RE1}
\end{equation}
for $k \neq 0$.

Finally, note that $\CC: \tilde{X}_\mathrm{r} \to \tilde{X}_\mathrm{r}$ is bounded, whence
$$\|\CC(\BB-\i k I)^{-1}\|_{\tilde{X}_\mathrm{r} \to \tilde{X}_\mathrm{r}} = O(|k|^{-1})$$
as $|k| \rightarrow \infty$.
Consequently $I-\CC(\BB-\i k I)^{-1}: \tilde{X}_\mathrm{r} \to \tilde{X}_\mathrm{r}$ is invertible for sufficiently large values of $|k|$ with
\begin{equation}
\|(I-\CC(\BB-\i k I)^{-1})^{-1}\|_{\tilde{X}_\mathrm{r} \to \tilde{X}_\mathrm{r}} = O(1) \label{RE2}
\end{equation}
as $|k| \rightarrow \infty$,
and the stated result follows from the identity
$$(L-\i k I)^{-1} = (\BB-\i k I)^{-1} \big(I-\CC(\BB-\i k I)^{-1}\big)^{-1}$$
and the estimates \eqn{RE1}, \eqn{RE2}.\qed

\begin{lemma} \label{Iooss condition}
The equation
\begin{equation}
Lw = -N(w^\dag) \label{ICAO}
\end{equation}
has a unique solution $w \in \tilde{D}_\star$ for each $w^\dag \in \tilde{D}_\star$ and the formula $w^\dag \mapsto w$ defines a smooth
mapping $\tilde{D}_\star \rightarrow \tilde{D}_\star$.
\end{lemma}
{\bf Proof.} Equation \eqn{ICAO} is equivalent to the equations
$$
u_1-u_{1xx} -6\sech^2 (x)u_1 = f(w^\dag),
$$
where $f(w^\dag)=(3\gamma_1+\gamma_2) \sech(x)(u_1^\dag)^2 +\gamma_2 u_1^\dag\phi_x^\dag +\gamma_1 (u_1^\dag)^3$, and
$$\phi_x = (u_1^\dag)^2+2\sech (x)u_1, \qquad v_1=0, \qquad \psi=0.$$
The result thus follows from Proposition \ref{Drazin prop}(i) and the fact that
$f$ and $(u_1,u_1^\dag) \mapsto (u_1^\dag)^2+2\sech (x)u_1$ define smooth
mappings $\tilde{D}_\star \to L^2_\mathrm{e}(\R)$ and $H^1_\mathrm{e}(\R) \times H^1_\mathrm{e}(\R)
\to H^1_\mathrm{e}(\R)$.\qed

\section{Transverse linear instability} \label{Apply GH}

Finally, we demonstrate the transverse linear instability of the line soliton using the
following general result due to Godey \cite{Godey16}.

\begin{theorem}[Godey]
Consider the differential equation
\begin{equation}
v_\tau = T v_t  + Lv , \label{Linearised multi-DE}
\end{equation}
in which $v(\tau,t)$ belongs to a real Banach space $\XX$. Suppose that $\YY$, $\ZZ$ are further real Banach spaces with the properties that
\begin{list}{(\roman{count})}{\usecounter{count}}
\item
$L: \ZZ \subseteq \XX \rightarrow \XX$ and $T:\YY \subseteq \XX \rightarrow \XX$ are closed linear operators with $\ZZ \subseteq \YY$,
\item
the equation is reversible: there exists an involution $S \in \LL(\XX) \cap \LL(\YY) \cap \LL(\ZZ)$ with $LSv=-SLv$ and
$TSv=-STv$ for all $v \in \ZZ$,
\item
$L$ has a pair $\pm \i \omega_0$ of isolated purely imaginary eigenvalues with odd algebraic multiplicity.
\end{list}
Under these hypotheses equation \eqn{Linearised multi-DE}
has a solution of the form $\e^{\lambda t} v_\lambda(\tau)$, where
$v_\lambda \in C^1({\mathbb R}, \XX) \cap C(\R, \ZZ)$ is periodic,  for each sufficiently small positive
value of $\lambda$; its period tends to $2\pi/\omega_0$ as $\lambda \rightarrow 0$.
\end{theorem}

Theorem \ref{Main result 2} is proved by applying Godey's theorem to \eqn{TS}, taking
$\XX=X$, $\YY=X$ and $\ZZ=D$ (and of course $\tau=y$ and $S(u_1,v_1,u_2,v_2,\phi,\psi)=(u_1,-v_1,u_2,-v_2,\phi,-\psi)$). The spectral hypothesis (iii) follows from Corollary \ref{Imaginary result}.
\\[2mm]
{\bf Acknowledgements}
E.~W.~was supported by the Swedish Research Council (grant no.~621-2012-3753).

\end{document}